\documentclass{iopart}
\usepackage{amssymb,mathdots}
\usepackage{iopams}
\usepackage{amsthm}
\usepackage{amscd}
\usepackage{amsfonts}
\usepackage{dsfont}
\usepackage{cite}
\usepackage{url}
\usepackage{graphicx}
\usepackage{mathrsfs}
\usepackage{pstricks,pst-plot,pst-node,pst-text,pst-3d}
\newtheorem{theorem}{Theorem}[section]

\newtheorem{definition}[theorem]{Definition}
\newtheorem{example}[theorem]{Example}
\newtheorem{remark}[theorem]{Remark}

\def\zede{\textit{\texttt{z}}}

\begin{document}
\title{Moving frames and conservation laws for Euclidean invariant Lagrangians}
\author{T\^ania M N Gon\c calves  and Elizabeth L Mansfield}
\address{School of Mathematics, Statistics and Actuarial Science, \\
University of Kent, Canterbury, CT2 7NZ, UK}
\eads{\mailto{T.M.N.Goncalves@kent.ac.uk},\quad \mailto{E.L.Mansfield@kent.ac.uk}}

\begin{abstract}
Noether's First Theorem yields conservation laws for Lagrangians with a variational symmetry group. The explicit formulae for the laws are well known and the symmetry group is known to act on the linear space generated by the conservation laws. In recent work the authors showed the mathematical structure behind both the Euler-Lagrange system and the set of conservation laws, in terms of the differential invariants of the group action and a moving frame. In this paper we demonstrate that the knowledge of this structure considerably eases finding the extremal curves for variational problems invariant under the special Euclidean groups $SE(2)$ and $SE(3)$.
\end{abstract}\vspace{-3ex}
 \submitto{\JPA} 

\section{Introduction}
In 1918 Emmy Noether wrote the seminal paper ``Invariante Variationsprobleme'' \cite{Noether}, where she showed that for differential systems derived from a variational principle, conservation laws could be obtained from Lie group actions which left the functional invariant.

Recently in \cite{GoncalvesMansfield}, it was proved that for Lagrangians that are invariant under some Lie symmetry group, Noether's conservation laws can be written in terms of a moving frame and vectors of invariants. Furthermore, in \cite{GoncalvesMansfield} the authors showed that for one-dimensional Lagrangians that are invariant under a semisimple Lie symmetry group, the new format for the conservation laws could reduce considerably the calculations needed to solve for the extremal curves. In particular, a classification is given for variational problems which are invariant under the three inequivalent $SL(2,\mathbb{C})$ actions on the plane (classified by Lie \cite{Lie}). In this paper we use the new structure of Noether's conservation laws presented in \cite{GoncalvesMansfield} to simplify variational problems that are invariant under $SE(2)$ and $SE(3)$; these groups are not semisimple.

In Section \ref{SNoetherCL}, we will briefly give an overview on moving frames, on differential invariants of a group action and on invariant calculus of variations. Throughout Section \ref{SNoetherCL} we will use the group action of $SE(2)$ on the plane as our pedagogical example.

In Section \ref{SE2}, we show in some detail how to compute the new version of Noether's conservation laws presented in Theorem \ref{tmel} for one-dimensional variational problems that are invariant under $SE(2)$, and then demonstrate how their invariantized Euler-Lagrange equations and Noether's conservation laws can be used to solve the integration problem.

Finally in Section \ref{SE3}, we present the simplified solution to the physically important one-dimensional variational problems that are left unchanged under the $SE(3)$ group action. 

\section{Structure of Noether's Conservation Laws}\label{SNoetherCL}
In this section, we will give a brief overview of concepts regarding moving frames, differential invariants of a group action and the invariant calculus of variations needed to understand the statements of our result. For more information on these subjects, see Fels and Olver \cite{FelsOlver,FelsOlverII}, Mansfield \cite{Mansfield} and Kogan and Olver \cite{KoganOlver}. We will use the $SE(2)$ action on the plane as our pedagogical example.

\subsection{Moving frames and differential invariants of a group action}
Here we are using moving frames as reformulated by Fels and Olver \cite{FelsOlver,FelsOlverII}, adapted to the context of differential algebra. 

Let $X$ be the space of independent variables with coordinates $\mathbf{x}=(x_1,...,x_p)$ and $U$ the space of dependent variables with coordinates $\mathbf{u}=(u^1,...,u^q)$. We will use a multiindex notation to represent the derivatives of $u^\alpha$, e.g. 
\begin{eqnarray}
u^\alpha_K=\frac{\partial^{|K|} u^\alpha}{\partial x_1^{k_1}\partial x_2^{k_2}\cdots\partial x_p^{k_p}},\nonumber
\end{eqnarray}
where the tuple $K=(k_1,...,k_p)$, represents a multiindex of differentiation of order $|K|=k_1+k_2+\cdots +k_p$. Hence, let $M=J^n(X\times U)$ be the $n$-th jet bundle with coordinates
\begin{eqnarray}
\zede=(x_1,...,x_p,u^1,...,u^q,u^1_1,...).\nonumber
\end{eqnarray}
On this space, the operator $\partial/\partial x_i$ extends to the \emph{total differentiation operator}
\begin{eqnarray}
D_i=\frac{D}{Dx_i}=\frac{\partial}{\partial x_i}+\sum_{\alpha=1}^q\sum_Ku^\alpha_{Ki}\frac{\partial}{\partial u^\alpha_K}.\nonumber
\end{eqnarray}

In this paper we are interested in using Noether's conservation laws to find the solutions that extremize variational problems which are invariant under $SE(2)$ and $SE(3)$. Thus, consider a Lagrangian $L(\zede)\mathrm{d}\mathbf{x}$ that is invariant under some symmetry group. Let a group $G$ act on the space $M$ as follows
\begin{eqnarray}
\begin{array}{ccl}
G\times M & \rightarrow & M\\
\zede & \mapsto & \widetilde{\zede}=g\cdot \zede,
\end{array}\nonumber
\end{eqnarray}
which satisfies either $g\cdot (h\cdot \zede)=(gh)\cdot \zede$, called a \emph{left action}, or $g\cdot (h\cdot \zede)=(hg)\cdot \zede$, called a \emph{right action}. We say a Lagrangian $L(\zede)\mathrm{d}\mathbf{x}$ is \emph{invariant} under some group action if
\begin{eqnarray}
L(\zede)\mathrm{d}\mathbf{x}=L(\widetilde{\zede})\mathrm{d}\widetilde{\mathbf{x}}\nonumber
\end{eqnarray}
for all $g \in G$.

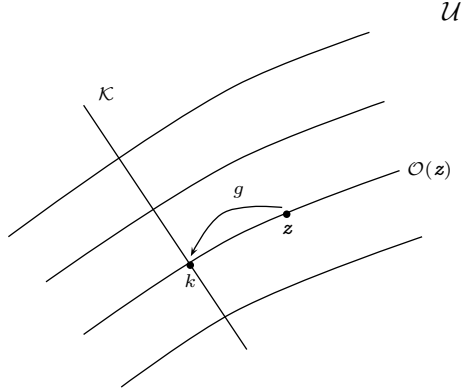
\begin{figure*}[h]
\begin{center}
\begin{pspicture}(-0.5,0)(5,5)
\pscurve[linewidth=0.5pt](1,0)(2.5,1)(5,2)
\pscurve[linewidth=0.5pt](0.5,0.7)(2.5,2)(4.7,2.88)
\pscurve[linewidth=0.5pt](0,1.4)(2.5,3)(4.5,3.8)
\pscurve[linewidth=0.5pt](-0.5,2)(2.5,4)(4.3,4.72)
\psline[linewidth=0.5pt](2.67,0.5)(0.5,3.75)
\pscircle*(3.2,2.3){0.05}
\pscircle*(1.918,1.622){0.05}
\pscurve[linewidth=0.5pt,arrowsize=2.5pt]{->}(3.14,2.39)(2.735,2.4)(2.33,2.3)(1.925,1.74)
\rput(3.2,2.1){$\scriptscriptstyle{\footnotesize{\zede}}$}
\rput(1.918,1.422){$\scriptstyle{k}$}
\rput(2.559,2.6){$\scriptstyle{g}$}
\rput(0.8,3.9){$\scriptstyle{\mathcal{K}}$}
\rput(5.1,2.88){$\scriptstyle{\mathcal{O}(\footnotesize{\zede})}$}
\rput(5.4,5){$\mathcal{U}$}
\end{pspicture}
\caption{A local foliation with a transverse cross section}\label{cartanmoving}
\end{center}
\end{figure*}

Consider a Lie group $G$ acting smootlhy on $M$ such that the action is \emph{free} and \emph{regular}. Then for every $\zede \in M$ there exists a neighbourhood $\mathcal{U}$ of $\zede$, as illustrated in Figure \ref{cartanmoving}, such that
\begin{itemize}
\item[-]
the group orbits have the dimension of the group $G$ and folliate $\mathcal{U}$;
\item[-]
there is a surface $\mathcal{K}\subset \mathcal{U}$ which intersects the group orbits transversally at a single point. This surface is called the \emph{cross section};
\item[-]
if $\mathcal{O}(\zede)$ represents the group orbit through $\zede$, then the group element $g\in G$ taking $\zede\in \mathcal{U}$ to $k$ is unique.
\end{itemize} 

Under these conditions, we can define a \emph{right moving frame} as the map $\rho: \mathcal{U}\rightarrow G$ which sends $\zede \in \mathcal{U}$ to the unique element $\rho(\zede)\in G$ such that
\begin{eqnarray}
\rho(\zede)\cdot \zede=k,\qquad \{k\}=\mathcal{O}(\zede)\cap\mathcal{K}.\nonumber
\end{eqnarray}
The element $g\in G$ in Figure \ref{cartanmoving} corresponds to $\rho(\zede)$.

To obtain the right moving frame, which sends $\zede$ to $k$, we must first define the cross section $\mathcal{K}$ as the locus of the set of equations $\psi_j(\zede)=0$, $j=1,...,r$, where $r$ is the dimension of $G$. Normally, the cross section is chosen so as to ease the calculations. Then solving the set of equations
\begin{equation}\label{normaeqns}
\psi_j(\widetilde{\zede})=\psi_j(g\cdot \zede)=0,\qquad j=1,...,r,
\end{equation}
known as the \emph{normalization equations}, for the $r$ group parameters describing $G$ yields the right moving frame in parametric form. Hence, the frame obtained satisfies 
\begin{eqnarray}
\psi_j(\rho(\zede)\cdot\zede)=0,\qquad j=1,...,r.\nonumber
\end{eqnarray}

By the implicit function theorem, a unique solution of (\ref{normaeqns}) provides an equivariant map, i.e. for a left action
\begin{eqnarray}
\rho(\widetilde{\zede})=\rho(\zede)g^{-1}\nonumber
\end{eqnarray}
and for a right action 
\begin{eqnarray}
\rho(\widetilde{\zede})=g^{-1}\rho(\zede).\nonumber
\end{eqnarray}

\begin{example}\label{example1}
Consider the $SE(2)$ group acting on curves in the $(x,y(x))$-plane as follows,
\begin{equation}\label{se2grpact}
\left(\begin{array}{c}x\\y\end{array}\right)\mapsto \left(\begin{array}{c}\widetilde{x}\\\widetilde{y}\end{array}\right)=\left(\begin{array}{cc}\cos{\theta}&\sin{\theta}\\ -\sin{\theta}&\cos{\theta}\end{array}\right)\left(\begin{array}{c}x-a\\y-b\end{array}\right),
\end{equation}
where $\theta$, $a$ and $b$ are constants that parametrize the group action. Here we are using the inverse action because it simplifies the calculations.

There is an induced action on the derivatives $y_K$, where $K$ is the index of differentiation with respect to $x$, called the \emph{prolonged action}. The induced action on $y_x$ is defined to be 
\begin{eqnarray}
\widetilde{y_x}=g\cdot y_x=\frac{\mathrm{d}\widetilde{y}}{\mathrm{d}\widetilde{x}}=\frac{\mathrm{d}\widetilde{y}/\mathrm{d}x}{\mathrm{d}\widetilde{x}/\mathrm{d}x}\nonumber
\end{eqnarray}
by the chain rule, so the action of (\ref{se2grpact}) on $y_x$ is 
\begin{eqnarray}
\widetilde{y_x}=\frac{-\sin{\theta}+y_x\cos{\theta}}{\cos{\theta}+y_x\sin{\theta}}.\nonumber
\end{eqnarray}
Similarly,
\begin{equation}\label{2ndpro}
\widetilde{y_{xx}}=g\cdot y_{xx}=\frac{\mathrm{d}^2\widetilde{y}}{\mathrm{d}(\widetilde{x})^2}=\frac{1}{\mathrm{d}\widetilde{x}/\mathrm{d}x}\frac{\mathrm{d}}{\mathrm{d}x}\left(\frac{\mathrm{d}\widetilde{y}}{\mathrm{d}\widetilde{x}}\right)=\frac{y_{xx}}{(\cos{\theta}+y_x\sin{\theta})^3}.
\end{equation}
If we consider $M$ to be the space with coordinates $(x,y,y_x,y_{xx},...)$, then the action is locally free near the identity of $SE(2)$. Thus, taking the normalization equations to be $\widetilde{x}=0$, $\widetilde{y}=0$ and $\widetilde{y_x}=0$, we obtain
\begin{equation}\label{se2frame}
a=x,\qquad b=y,\qquad and \qquad \theta=\arctan{y_x}
\end{equation}
as the frame in parametric form.
\end{example}

\begin{remark}
In this paper we will consider all independent variables to be invariant. If these are not invariant, then we can reparametrize and set the original independent variables as depending on the new invariant parameters.
\end{remark}

\begin{theorem}
Let $\rho(\zede)$ be a right moving frame. Then the quantity $I(\zede)=\rho(\zede)\cdot\zede$ is an invariant of the group action (see \cite{FelsOlver}).
\end{theorem}
Consider $\zede=(z_1,...,z_m) \in M$ and let the normalization equations $\widetilde{z_i}=c_i$ for $i=1,...,r$, where $r$ is the dimension of the group $G$, then
\begin{eqnarray}
\rho(\zede)\cdot\zede=(c_1,...,c_r,I(\zede_{r+1}),...,I(\zede_m)),\nonumber
\end{eqnarray}
where 
\begin{eqnarray}
I(\zede_l)=g\cdot\zede|_{g=\rho(\footnotesize{\zede})},\qquad l=r+1,...,m.\nonumber
\end{eqnarray}

\vspace{0.1cm}
\noindent\textbf{Example 2.1 (cont.)} \textit{Evaluating $\widetilde{y_{xx}}$ given by (\ref{2ndpro}) at the frame (\ref{se2frame}) yields 
\begin{eqnarray}
\widetilde{y_{xx}}|_{(a=x,b=y,\theta=\arctan{y_x})}=\frac{y_{xx}}{(1+y_x^2)^{3/2}},\nonumber
\end{eqnarray}
the \emph{Euclidean curvature}. So evaluating $\widetilde{y_K}$ at the frame (\ref{se2frame}) yields a differential invariant.}

\begin{definition}
For any prolonged action in the jet space $J^n(X\times U)$, the invariantized jet coordinates are denoted as
\begin{equation}\label{formulainv}
J_i=I(x_i)=\widetilde{x_i}|_{g=\rho(\footnotesize{\zede})},\qquad I^\alpha_K=I(u^\alpha_K)=\widetilde{u^\alpha_K}|_{g=\rho(\footnotesize{\zede})}.
\end{equation}
These are also known as the \emph{normalized differential invariants}.
\end{definition}

\begin{example}\label{example2}
Consider the group action of $SE(2)$ as in Example \ref{example1}. Since $x$ is not invariant we reparametrize $(x,y(x))$ as $(x(s),y(s))$, where $s$ is invariant and let $g\in SE(2)$ act on $(x(s),y(s))$ as in Example \ref{example1}. Solving the normalization equations $\widetilde{x}=0$, $\widetilde{y}=0$ and $\widetilde{y_s}=0$, we obtain the frame 
\begin{equation}\label{frameSE2}
a=x,\qquad b=y,\qquad \theta=\arctan{\left(\frac{y_s}{x_s}\right)}.
\end{equation}
We have then
\begin{eqnarray}
g\cdot\zede|_{g=\rho(\footnotesize{\zede})}&=(\widetilde{s},\widetilde{x},\widetilde{y},\widetilde{x_s},\widetilde{y_s},\widetilde{y_{ss}})|_{g=\rho(\footnotesize{\zede})}\nonumber\\[10pt]
&=(I(s),I^x,I^y,I^x_1,I^y_1,I^y_{11})\nonumber\\\label{invs1}
&=\left(s,0,0,\sqrt{x_s^2+y_s^2},0,\frac{x_sy_{ss}-y_sx_{ss}}{\sqrt{x_s^2+y_s^2}}\right).
\end{eqnarray}
The second, third and fifth components of (\ref{invs1}) correspond to the normalization equations $\widetilde{x}=0$, $\widetilde{y}=0$ and $\widetilde{y_s}=0$ respectively. The fourth and sixth components, $I^x_1$ and $I^y_{11}$ respectively, are the lowest order differential invariants and all higher order invariants can be obtained in terms of them and their derivatives.
\end{example}

\begin{theorem}
(\textbf{Replacement Theorem} \cite{FelsOlverII}) If $f(\zede)$ is an invariant, then 
\begin{eqnarray}
f(\zede)=f(I(\zede)).\nonumber
\end{eqnarray}
\end{theorem}

This theorem allows one to find the $I^\alpha_K$ in terms of historically well-known invariants without having to solve for the frame.

\vspace{0.15cm}
\noindent \textbf{Example \ref{example2} (cont.)} \textit{We know that $SE(2)$ preserves $|\mathbf{x}_s|$, thus applying the Replacement Theorem we obtain}
\begin{eqnarray}
|\mathbf{x}_s|=\sqrt{x_s^2+y_s^2}=\sqrt{(I^x_1)^2+(I^y_1)^2}=|I^x_1|,\nonumber
\end{eqnarray}
\textit{which yields that $\sqrt{x_s^2+y_s^2}=I^x_1$ up to a sign. Next we know that the Euclidean curvature $\kappa$ is also invariant under $SE(2)$, and we obtain}
\begin{eqnarray}
\kappa=\frac{x_sy_{ss}-y_sx_{ss}}{(x_s^2+y_s^2)^{3/2}}=\frac{I^x_1I^y_{11}-I^y_1I^x_{11}}{((I^x_1)^2+(I^y_1)^2)^{3/2}}=\frac{I^y_{11}}{(I^x_1)^2},\nonumber
\end{eqnarray}
\textit{which gives us $I^y_{11}$ in terms of $\kappa$ and $|\mathbf{x}_s|$.}

\vspace{0.15cm}
Since we are considering all independent variables to be invariant, all total differential operators will also be invariant. Hence, the \emph{invariantized differential operators} 
\begin{eqnarray}
\mathcal{D}_i=\widetilde{D_i}|_{g=\rho(\scriptscriptstyle{\zede})}=D_i.\nonumber
\end{eqnarray}

We know that
\begin{eqnarray}
\frac{\partial }{\partial x_i}u^\alpha_K=u^\alpha_{Ki},\nonumber
\end{eqnarray}
although the same cannot be said about its invariantized version and in general $\mathcal{D}_iI^\alpha_K\ne I^\alpha_{Ki}$; indeed we have that
\begin{eqnarray}\label{diffinvariant}
\mathcal{D}_iI^\alpha_K=I^\alpha_{Ki}+M^\alpha_{Ki},
\end{eqnarray}
where $M^\alpha_{Ki}$ is known as the \emph{correction term}. We will not go into its calculation since it would take us to far afield, but for more information on correction terms see \S 4.5 \cite{Mansfield}. In any case, software exists to calculate these correction terms \cite{AIDA}. Equation \ref{diffinvariant} shows that the processes of invariantization and differentiation do not commute. Considering two generating differential invariants $I^\alpha_J$ and $I^\alpha_L$ and letting $JK=LM$ so that $I^\alpha_{JK}=I^\alpha_{LM}$, then this implies that
\begin{equation}
\mathcal{D}_KI^\alpha_J-M^\alpha_{JK}=\mathcal{D}_MI^\alpha_L-M^\alpha_{LM}.
\end{equation}
These equations are called \emph{syzygies} or \emph{differential identities}. These will play a crucial role in the obtention of the invariantized Euler-Lagrange equations and Noether's conservation laws.

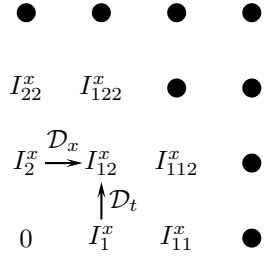
\begin{figure*}[h]
\begin{center}
\begin{pspicture}(,0)(4.5,3.5)
\rput(-1,0){$0$}
\rput(0,0){$I^x_1$}
\rput(1,0){$I^x_{11}$}
\rput(-1,1){$I^x_2$}
\rput(0,1){$I^x_{12}$}
\rput(0,2){$I^x_{122}$}
\rput(-1,2){$I^x_{22}$}
\rput(1,1){$I^x_{112}$}
\psline[arrowlength=2]{->}(-0.75,1)(-0.25,1)
\rput(-0.5,1.3){\textcolor{black}{$\mathcal{D}_x$}}
\psdots[linecolor=black,dotstyle=*,dotscale=2]
(-1,3)(0,3)(1,3)(2,3)(1,2)(2,2)(2,1)(2,0)
\psline[arrowlength=2]{->}(0,0.25)(0,0.75)
\rput(0.3,0.5){$\mathcal{D}_t$}
\end{pspicture}
\caption{Paths to $I^x_{12}$}\label{twopaths}
\end{center}
\end{figure*}

\vspace{0.15cm}
\noindent \textbf{Example \ref{example2} (cont.)} \textit{If we set $x=x(s,t)$ and $y=y(s,t)$ and take the normalization equations as before, we obtain}
\begin{eqnarray}
\widetilde{x_t}|_{g=\rho(\footnotesize{\zede})}=I^x_2=\frac{x_sx_t+y_sy_t}{\sqrt{x_s^2+y_s^2}},\qquad\widetilde{y_t}|_{g=\rho(\footnotesize{\zede})}=I^y_2=\frac{x_sy_t-y_sx_t}{\sqrt{x_s^2+y_s^2}}.\nonumber
\end{eqnarray}
\noindent\textit{Furthermore, since both $s$ and $t$ are invariant, $\mathcal{D}_s$ and $\mathcal{D}_t$ commute. From Figure \ref{twopaths}, we can see that there are two ways in which we can obtain $I^x_{12}$} 

\vspace{0.15cm}
\noindent\textit{and since both ways must be equal, we get a syzygy between $I^x_2$ and $\eta=I^x_1$. The syzygy is}
\begin{equation}\label{syzygy1}
\mathcal{D}_t\eta=\mathcal{D}_sI^x_2-\kappa\eta I^y_2.
\end{equation}
\noindent\textit{Similarly, we have a syzygy between $I^y_2$ and $I^y_{11}$ and the syzygy is}
\begin{equation}\label{syzygy2}
\mathcal{D}_tI^y_{11}=\mathcal{D}_s^2I^y_2-\frac{\eta_s}{\eta}\mathcal{D}_sI^y_2-\kappa^2 \eta^2 I^y_2+2\kappa \eta\mathcal{D}_sI^x_2+\kappa_s \eta I^x_2.
\end{equation}

\subsection{Invariant calculus of variations}\label{ICV}
Kogan and Olver in \cite{KoganOlver} studied invariant calculus of variations from a geometric point of view, we instead do it from a differential algebra point of view. 

Assume Lagrangians to be smooth functions of $\mathbf{x}$, $\mathbf{u}$ and  finitely many derivatives of the $u^\alpha$ and denote these as $\mathscr{L}[\mathbf{u}]=\int L[\mathbf{u}]\mathrm{d}\mathbf{x}$. Furthermore, suppose these are invariant under some group action and let the $\kappa_j$, $j=1,...,N$, denote the generating differential invariants of that group action. Also, assume that the action leaves the independent variables invariant so that the Lagrangians can be rewritten as $\int L[\boldsymbol{\kappa}]\mathrm{d}\mathbf{x}$. 

To obtain the invariantized Euler-Lagrange equations we proceed in a similar way as for finding the Euler-Lagrange equations in the original variables $(\mathbf{x},\mathbf{u})$.

Recall that if $\mathbf{x}\mapsto (\mathbf{x},\mathbf{u})$ extremizes the functional $\mathscr{L}[\mathbf{u}]$, then for a small perturbation of $\mathbf{u}$,
\begin{eqnarray}
 0&=\left.\frac{\mathrm{d}}{\mathrm{d}\varepsilon}\right|_{\varepsilon=0}\mathscr{L}[\mathbf{u}+\varepsilon\mathbf{v}]\nonumber \\
&=\int \sum_{\alpha=1}^q\left[\mathsf{E}^\alpha(L)v^\alpha+\sum_i\frac{D}{Dx_i}\left(\frac{\partial L}{\partial u^\alpha_i}v^\alpha+\cdots\right)\right]\mathbf{d}\mathbf{x}\nonumber
\end{eqnarray}
after differentiation under the integral sign and integration by parts, where
\begin{eqnarray}
\mathsf{E}^\alpha=\sum_K(-1)^{|K|}\frac{D^{|K|}}{D x_1^{k_1}\cdots D x_p^{k_p}}\frac{\partial}{\partial u^\alpha_K}\nonumber
\end{eqnarray}
is the Euler operator with respect to the dependent variables $u^\alpha$. The boundary terms correspond to Noether's conservation laws and the variation $\mathbf{v}$ to the infinitesimals.

To obtain the invariantized Euler-Lagrange equations, we first introduce a dummy invariant independent variable $t$ and set the $u^\alpha=u^\alpha(\mathbf{x},t)$. The introduction of this new independent variable results in $q$ new invariants $I^\alpha_t=g\cdot u^\alpha_t|_{g=\rho(\footnotesize{\zede})}$ and a set of syzygies $\mathcal{D}_t\boldsymbol{\kappa}=\mathcal{H}I(\mathbf{u}_t)$ that is
\begin{equation}
\mathcal{D}_t\left(\begin{array}{c}\kappa_1\\\vdots\\\kappa_N\end{array}\right)=\mathcal{H}\left(\begin{array}{c}I^1_t\\\vdots\\ I^q_t\end{array}\right),
\end{equation}
where $\mathcal{H}$ is a $N\times q$ matrix of operators depending only on the $\mathcal{D}_i$, for $i=1,...,p$, the $\kappa_j$, for $j=1,...,N$, and their invariant derivatives. Since all independent variables are invariant, we have that all differential operators commute, specifically $[\mathcal{D}_i,\mathcal{D}_t]=0$, for all $i=1,...,p$.

\begin{remark}
Up to this moment we have represented the independent variables as $\mathbf{x}$ and the dependent variables as $\mathbf{u}$. Since the examples in this paper only involve two independent variables, $s$ and $t$, we will represent the coordinates of the space of dependent variables as $\mathbf{x}$, where the dependent variables will be the usual space coordinates.
\end{remark}

For simplicity, consider a one-dimensional Lagrangian $L(x,y,y_x,y_{xx},...)\mathrm{d}x$ with a finite number of arguments which is invariant under the $SE(2)$ group action (\ref{se2grpact}). Such variational problems can be rewritten in terms of the generating invariants of its group action, in this case the Euclidean curvature, $\kappa$, and its derivatives with respect to $s$, the Euclidean arc length. We reparametrize $(x,y(x))$ as $(x(s),y(s))$, and to fix parametrization as arc length, we introduce $\eta=\sqrt{x_s^2+y_s^2}=1$ as a constraint. Thus, we consider the invariantized variational problem
\begin{eqnarray}\label{functionalse2}
\int [L(\kappa,\kappa_s,\kappa_{ss},...)-\lambda(s)(\eta-1)]\mathrm{d}s,
\end{eqnarray}
where $\lambda(s)$ is a Lagrange multiplier. This constraint does not reduce the solution set and it will simplify the calculations. Symbolically, we know that
\begin{eqnarray}
\left.\frac{\mathrm{d}}{\mathrm{d}\varepsilon}\right|_{\varepsilon=0}\mathscr{L}[\mathbf{x}+\varepsilon\mathbf{v}]=\left.\frac{D}{Dt}\right|_{\mathbf{x}_t=\mathbf{v}}\mathscr{L}[\mathbf{x}].\nonumber
\end{eqnarray}
Hence, after differentiating (\ref{functionalse2}) under the integral sign and integrating by parts we obtain
\begin{eqnarray}
\fl \mathcal{D}_t\int [L(\kappa,\kappa_s,\kappa_{ss},...)-\lambda(s)(\eta-1)]\mathrm{d}s\nonumber\\[5pt]
\fl\quad=\int \left[\frac{\partial L}{\partial \kappa}\mathcal{D}_t\kappa+\frac{\partial L}{\partial \kappa_s}\mathcal{D}_s\mathcal{D}_t\kappa+\frac{\partial L}{\partial \kappa_{ss}}\mathcal{D}_s^2\mathcal{D}_t\kappa+\cdots-\lambda(s)\mathcal{D}_t\eta\right]\mathrm{d}s\nonumber\\[5pt]
\fl \quad=\int \Big[\left(\frac{\partial L}{\partial \kappa}-\mathcal{D}_s\left(\frac{\partial L}{\partial \kappa_s}\right)+\mathcal{D}_s^2\left(\frac{\partial L}{\partial \kappa_{ss}}\right)-\mathcal{D}_s^3\left(\frac{\partial L}{\partial \kappa_{sss}}\right)+\cdots\right)\mathcal{D}_t\kappa-\lambda(s)\mathcal{D}_t\eta\nonumber\\[5pt]
\fl \quad\quad+\mathcal{D}_s\Big(\sum_{m=1}\sum_{n=0}^{m-1}(-1)^n\mathcal{D}_s^n\left(\frac{\partial L}{\partial \kappa_m}\right)\mathcal{D}_s^{m-1-n}\mathcal{D}_t\kappa\Big)\Big]\mathrm{d}s\nonumber\\[5pt]
\fl \quad=\int \Big[\mathsf{E}^\kappa(L)\pscirclebox[linewidth=0.3pt]{\mathcal{D}_t\kappa}-\lambda(s)\pscirclebox[linewidth=0.3pt]{\mathcal{D}_t\eta}+\mathcal{D}_s\Big(\sum_{m=1}\sum_{n=0}^{m-1}(-1)^n\mathcal{D}_s^n\left(\frac{\partial L}{\partial \kappa_m}\right)\mathcal{D}_s^{m-1-n}\mathcal{D}_t\kappa\Big)\Big]\mathrm{d}s,\nonumber
\end{eqnarray}
where 
\begin{eqnarray}
\kappa_m=\frac{\partial^m \kappa}{\partial s^m}.\nonumber
\end{eqnarray} 
Next we substitute the circled $\mathcal{D}_t\eta$ and $\mathcal{D}_t\kappa$ by their respective syzygies
\begin{equation}
\mathcal{D}_t\left(\begin{array}{c}\eta\\\kappa\end{array}\right)=\left(\begin{array}{cc}\mathcal{D}_s & -\kappa\\[10pt]
\kappa_s & \mathcal{D}_s^2+\kappa^2\end{array}\right)\left(\begin{array}{c}I^x_2\\[10pt]I^y_2\end{array}\right)=\left(\begin{array}{c}\mathcal{H}_1\\[10pt]\mathcal{H}_2\end{array}\right)I(\mathbf{x}_t),
\end{equation}
where we have already set $\eta=1$; it makes no difference to the final result to do this at this point. Hence,
\begin{eqnarray}
\fl \int \Big[\left(\mathsf{E}^\kappa(L)\mathcal{H}_2-\lambda(s)\mathcal{H}_1\right)I(\mathbf{x}_t)+\mathcal{D}_s\Big(\sum_{m=1}\sum_{n=0}^{m-1}(-1)^n\mathcal{D}_s^n\left(\frac{\partial L}{\partial \kappa_m}\right)\mathcal{D}_s^{m-1-n}\mathcal{D}_t\kappa\Big)\Big]\mathrm{d}s\nonumber\\
\fl \quad=\int\Big[\left(\mathcal{H}_2^\ast(\mathsf{E}^\kappa(L))-\mathcal{H}_1^\ast(\lambda(s))\right)I(\mathbf{x}_t)+\mathcal{D}_s\Big(-\lambda(s)I^x_2+\mathsf{E}^\kappa(L)\mathcal{D}_sI^y_2-\mathcal{D}_s\mathsf{E}^\kappa(L)I^y_2\nonumber\\\label{boundarytermsse2}
\fl \quad\quad+\sum_{m=1}\sum_{n=0}^{m-1}(-1)^n\mathcal{D}_s^n\left(\frac{\partial L}{\partial \kappa_m}\right)\mathcal{D}_s^{m-1-n}\mathcal{D}_t\kappa\Big)\Big]\mathrm{d}s,
\end{eqnarray}
after a second set of integration by parts, and where $\mathcal{H}_1^\ast$ and $\mathcal{H}_2^\ast$ are the respective adjoint operators of $\mathcal{H}_1$ and $\mathcal{H}_2$. The vector $I(\mathbf{x}_t)$ corresponds to the variation $\mathbf{v}$ and thus the coefficients of $I^x_2$ and $I^y_2$ represent respectively the Euler-Lagrange equations
\begin{eqnarray}
\mathsf{E}^x(L)=\kappa_s\mathsf{E}^\kappa(L)+\lambda_s,\\
\mathsf{E}^y(L)=\mathcal{D}_s^2\mathsf{E}^\kappa(L)+\kappa^2\mathsf{E}^\kappa(L)+\lambda(s)\kappa.
\end{eqnarray}
 
We will use $\mathsf{E}^x(L)=0$ to eliminate $\lambda(s)$ instead of $\mathsf{E}^y(L)=0$, as it contains derivatives of $\kappa$ of lower order. Since $L$ does not depend on $s$ explicitly, by the result in page 220 of \cite{Mansfield}, $\kappa_s\mathsf{E}^\kappa(L)$ is a total derivative, specifically
\begin{eqnarray}
\kappa_s\mathsf{E}^\kappa(L)=\mathcal{D}_s\left(L-\sum_{m=1}\sum_{j=0}^{m-1}(-1)^{j}\mathcal{D}_s^{j}\left(\frac{\partial L}{\partial \kappa_m}\right)\kappa_{m-j}\right),\nonumber
\end{eqnarray}
then we obtain
\begin{equation}\label{constraint1}
\lambda(s)=-L+\sum_{m=1}\sum_{j=0}^{m-1}(-1)^{j}\mathcal{D}_s^{j}\left(\frac{\partial L}{\partial \kappa_m}\right)\kappa_{m-j},
\end{equation}
where the constant of integration has been absorbed into $\lambda(s)$ (see Remark 7.1.9. of \cite{Mansfield}). Hence, we are left with one invariantized Euler-Lagrange equation in one unknown
\begin{equation}
\fl\mathsf{E}^y(L)=\mathcal{D}_s^2\mathsf{E}^\kappa(L)+\kappa^2\mathsf{E}^\kappa(L)-\kappa\left(L-\sum_{m=1}\sum_{j=0}^{m-1}(-1)^{m}\mathcal{D}_s^{m}\left(\frac{\partial L}{\partial \kappa_m}\right)\kappa_{m-j}\right).
\end{equation}
This particular equation agrees with the one appearing in Kogan and Olver \cite{KoganOlver}.

\subsection{New version of Noether's conservation laws}\label{NVNcl}
If one calculates the conservation laws for one-dimensional variational problems
that are invariant under $SE(2)$ from Noether's First Theorem (for the formulae of these, see Theorem 4.29 of \cite{Olver}) and then rewrite these in terms of the Euclidean curvature $\kappa$ and its derivatives with respect to the Euclidean arc length $s$, one obtains
\begin{equation}\label{lawsSE2}
\fl\underbrace{\left(\begin{array}{ccc}\cos{\theta} &-\sin{\theta} & 0\\\sin{\theta} &\cos{\theta} & 0\\a\sin{\theta}-b\cos{\theta}&a\cos{\theta}+b\sin{\theta}&1\end{array}\right)}_{R(g)}\Bigg|_{g=\rho(\footnotesize{\zede})^{-1}}\underbrace{\left(\begin{array}{c}-\lambda(s)-\kappa\mathsf{E}^\kappa(L)\\-\mathcal{D}_s\mathsf{E}^\kappa(L)\\
\mathsf{E}^\kappa(L)\end{array}\right)}_{\boldsymbol{\upsilon}(I)}=\mathbf{c},
\end{equation} 
where $\rho(\zede)^{-1}$ is the right moving frame (\ref{frameSE2}), $\lambda(s)$ is the Lagrange multiplier obtained in (\ref{constraint1}) and $\mathbf{c}$ the constant vector. 

We note that $R(g)$ is the Adjoint representation of $SE(2)$ on its Lie algebra $\frak{se}(2)$. To see how calculations in Section \ref{ICV} can yield the result in (\ref{lawsSE2}), we first show how the Adjoint representation is calculated in the context we will need.

Consider the $SE(2)$ group action 
\begin{equation}\label{leftse2grpact}
\left(\begin{array}{c}\widetilde{x}\\\widetilde{y}\end{array}\right)=\left(\begin{array}{cc}\cos{\theta} & -\sin{\theta}\\
\sin{\theta} & \cos{\theta}\end{array}\right)\left(\begin{array}{c}x\\y\end{array}\right)+\left(\begin{array}{c}a\\b\end{array}\right)
\end{equation}
with generating infinitesimal vector fields
\begin{eqnarray}
\partial_x,\qquad \partial_y,\qquad -y\partial_x+x\partial_y.\nonumber
\end{eqnarray}
Let $g\in SE(2)$ act on 
\begin{eqnarray}
\mathbf{v}=\alpha\partial_x+\beta\partial_y+\gamma(-y\partial_x+x\partial_y)\nonumber
\end{eqnarray}
as in (\ref{leftse2grpact}), where $\alpha$, $\beta$ and $\gamma$ are constants. Hence,
\begin{eqnarray}
\fl \begin{array}{rl} g\cdot \mathbf{v}&\kern-8pt=\alpha\partial_{\widetilde{x}}+\beta\partial_{\widetilde{y}}+\gamma(-\widetilde{y}\partial_{\widetilde{x}}+\widetilde{x}\partial_{\widetilde{y}})\\[10pt]
&\kern-8pt=\left(\begin{array}{ccc}\alpha&\beta&\gamma\end{array}\right)\left(\begin{array}{ccc}\cos{\theta} &-\sin{\theta} & 0\\\sin{\theta} &\cos{\theta} & 0\\a\sin{\theta}-b\cos{\theta}&a\cos{\theta}+b\sin{\theta}&1\end{array}\right)\left(\begin{array}{c}\partial_x\\\partial_y\\-y\partial_x+x\partial_y\end{array}\right).\end{array}\nonumber
\end{eqnarray}
Thus, $R(g)$ is the \emph{Adjoint representation of} $SE(2)$, denoted as $\mathcal{A}d(g)$.

\begin{remark}
In Example \ref{example1} we used the right action of $SE(2)$ on the plane to calculate the right moving frame, as it simplified its calculation. However, to compute $\mathcal{A}d(\rho)^{-1}$ we considered the left action of $SE(2)$ on the plane, avoiding in this way the need to calculate the inverse of $\mathcal{A}d(\rho)$.
\end{remark}

Next recall the boundary terms (\ref{boundarytermsse2}) obtained in the calculation of the invariantized Euler-Lagrange equations,
\begin{eqnarray}
\fl -\lambda(s)I^x_2+\mathsf{E}^\kappa(L)\mathcal{D}_sI^y_2-\mathcal{D}_s\mathsf{E}^\kappa(L)I^y_2+\sum_{m=1}\sum_{n=0}^{m-1}(-1)^n\mathcal{D}_s^n\left(\frac{\partial L}{\partial \kappa_m}\right)\mathcal{D}_s^{m-1-n}\mathcal{D}_t\kappa=c,\nonumber
\end{eqnarray}
where $c$ is a constant.
Substituting $\mathcal{D}_sI^y_2$ and $\mathcal{D}_s^{m-1-n}\mathcal{D}_t\kappa$ for all $m$ in the above expression by the differential formulae 
\begin{eqnarray}
\mathcal{D}_sI^y_2&=I^y_{12}-\kappa I^x_2,\nonumber\\
\mathcal{D}_t\kappa&=I^y_{112}-2\kappa I^x_{12},\nonumber\\
\mathcal{D}_s\mathcal{D}_t\kappa&=I^y_{1112}-2\kappa^2 I^y_{12}-3\kappa I^x_{112}-2\kappa_s I^x_{12},\nonumber\\
&\qquad\vdots\nonumber
\end{eqnarray}
where these were obtained from (\ref{diffinvariant}), and rewriting it as 
\begin{equation}
\fl(\kern-5pt\begin{array}{ccc}I^x_2&\kern-8pt I^x_{12}&\kern-8pt\cdots\end{array}\kern-5pt)\underbrace{\left(\kern-8pt\begin{array}{c}-\lambda(s)-\kappa\mathsf{E}^\kappa(L)\\[4pt]
-2\kappa\frac{\partial L}{\partial \kappa_s}-2\kappa_s\frac{\partial L}{\partial \kappa_{ss}}+\cdots\\[4pt]
-3\kappa\frac{\partial L}{\partial \kappa_{ss}}+\cdots\\[4pt]
\vdots\end{array}\kern-8pt\right)}_{\mathcal{C}^x}+(\kern-5pt\begin{array}{ccc}I^y_2&\kern-8pt I^y_{12}&\kern-8pt\cdots\end{array}\kern-5pt)\underbrace{\left(\kern-8pt\begin{array}{c}-\mathcal{D}_s\mathsf{E}^\kappa(L)\\[4pt]\mathsf{E}^\kappa(L)-2\kappa^2\frac{\partial L}{\partial\kappa_{ss}}+\cdots\\[4pt]\frac{\partial L}{\partial \kappa_s}+\cdots\\[4pt]
\frac{\partial L}{\partial\kappa_{ss}}+\cdots\\[4pt]
\vdots\end{array}\kern-8pt\right)}_{\mathcal{C}^y}=k\nonumber
\end{equation}
yields the boundary terms in a form that is linear in the $I^\alpha_{2K}$. 

We now let $t$ be a group parameter. If the parameters are $(a_1,...,a_r)$ and $t=a_j$, then from Theorem 3 of \cite{GoncalvesMansfield}, it is shown that the vectors $(\begin{array}{ccc}I^\alpha_2& I^\alpha_{12}&\cdots\end{array})$ can be written as the product of row $j$ of $\mathcal{A}d(\rho)^{-1}$ and the matrix of invariantized infinitesimals
\begin{eqnarray}
\Omega^\alpha(I)=\left(\widetilde{\zeta}^i_j(I)\right),\qquad \zeta^i_j=\left.\frac{\partial \widetilde{z_i}}{\partial a_j}\right|_{g=e},\nonumber
\end{eqnarray}
where $\alpha$ represents a dependent variable, $a_j$ a group parameter, and $e$ the identity element. The vector of invariants in (\ref{lawsSE2}) equals the sum of the products of the matrices of invariantized infinitesimals $\Omega^\alpha(I)$ with the vectors $\mathcal{C}^\alpha$,
\begin{eqnarray}
\Omega^x(I)\mathcal{C}^x+\Omega^y(I)\mathcal{C}^y=\left(\begin{array}{c}-\lambda(s)-\kappa\mathsf{E}^\kappa(L)\\-\mathcal{D}_s\mathsf{E}^\kappa(L)\\
\mathsf{E}^\kappa(L)\end{array}\right),\nonumber
\end{eqnarray}
where
\begin{eqnarray}
\fl \Omega^x(I)=\bordermatrix{ & x & x_s & x_{ss} & \cdots \cr a & 1 & 0 & 0 & \cdots\cr b & 0 & 0 & 0 & \cdots \cr \theta & 0 & 0 & -\kappa & \cdots},\qquad \Omega^y(I)=\bordermatrix{ & y & y_s & y_{ss} & y_{sss} & \cdots \cr a & 0 & 0 & 0 & 0 & \cdots\cr b & 1 & 0 & 0 & 0 & \cdots \cr \theta & 0 & 1 & 0 & -\kappa^2 & \cdots}.\nonumber
\end{eqnarray}

Noether's conservation laws for one-dimensional Lagrangians invariant under $SE(2)$ can be written as
\begin{eqnarray}\label{conslawsse2}
\left(\begin{array}{ccc}
x_s & -y_s & 0\\
y_s & x_s & 0\\
xy_s-yx_s & xx_s+yy_s & 1
\end{array}\right)\left(\begin{array}{c}-\lambda(s)-\kappa\mathsf{E}^\kappa(L)\\-\mathcal{D}_s\mathsf{E}^\kappa(L)\\
\mathsf{E}^\kappa(L)\end{array}\right)=\mathbf{c},
\end{eqnarray}
where $x_s^2+y_s^2=1$ was used to simplify the conservation laws and
\begin{eqnarray}
\lambda(s)=-L+\sum_{m=1}\sum_{j=0}^{m-1}(-1)^{j}\mathcal{D}_s^{j}\left(\frac{\partial L}{\partial \kappa_m}\right)\kappa_{m-j}.\nonumber
\end{eqnarray}
The following theorem generalizes what we have just seen for one-dimensional variational problems that are invariant under $SE(2)$; it states that the conservation laws from Noether's First Theorem can be written as the divergence of the product of a moving frame with vectors of invariants.

\begin{theorem}\label{tmel}
Let $\int L(\kappa_1,\kappa_2,...)\mathrm{d}\mathbf{x}$ be invariant under $G\times M\rightarrow M$, where $M=J^n(X\times U)$, with generating invariants $\kappa_j$, for $j=1,...,N$, and let $\widetilde{x_i}=g\cdot x_i=x_i$, for $i=1,...,p$. Let $(a_1,...,a_r)$ be coordinates of $G$ near the identity $e$, and $\mathbf{v}_i$, for $i=1,...,r$, the associated infinitesimal vector fields. Furthermore, let $\mathcal{A}d(g)$ be the Adjoint representation of $G$ with respect to these vector fields. For each dependent variable, define the matrix of infinitesimals to be
\begin{eqnarray}
\Omega^\alpha(\widetilde{\zede})=\left(\widetilde{\zeta^i_j}\right),\nonumber
\end{eqnarray}
where
\begin{eqnarray}
\zeta^i_j=\left.\frac{\partial \widetilde{\zede_i}}{\partial a_j}\right|_{g=e}\nonumber
\end{eqnarray}
are the infinitesimals of the prolonged group action. Let $\Omega^\alpha(I)$, for $\alpha=1,...,q$, be the invariantized version of the above matrices.

Introduce a dummy invariant variable $t$ to effect the variation and then integration by parts yields
\begin{eqnarray}
\displaystyle{\frac{\partial}{\partial t}\int L(\kappa_1,\kappa_2,...)\mathrm{d}\mathbf{x}=\int\Big[\sum_\alpha\mathsf{E}^\alpha(L)I_t^\alpha+\mathsf{Div}(P)\Big]\mathrm{d}\mathbf{x}},\nonumber
\end{eqnarray}
where this defines the $p$-tuple $P$, whose components are of the form
\begin{eqnarray}
P_i=\displaystyle{\sum_{\alpha,J} I^\alpha_{t K}C^\alpha_{i,K},\qquad i=1,...,p,}\nonumber
\end{eqnarray}
and the vectors $\mathcal{C}^\alpha_{i}=(C^\alpha_{i,K})$. Let $\rho(\zede)^{-1}$ be a right frame with canonical invariants $I^\alpha_{t K}=I(u^\alpha_{t K})$, where $K$ is the index of differentiation with respect to the independent variables $x_i$, for $i=1,...,p$. Then the $r$ conservation laws obtained via Noether's First Theorem can be written in the form
\begin{equation}
\sum_i\frac{D}{Dx_i}\mathcal{A}d(\rho(\zede))^{-1}\boldsymbol{\upsilon}_i(I)=0,
\end{equation}
where
\begin{eqnarray}
\boldsymbol{\upsilon}_i(I)=\sum_\alpha\Omega^\alpha(I)\mathcal{C}^\alpha_i.\nonumber
\end{eqnarray}
\end{theorem}

The proof can be found in \cite{GoncalvesMansfield}.
\begin{remark}
One can notice that none of the summands of
\begin{eqnarray}
\sum_{m=1}\sum_{n=0}^{m-1}(-1)^n\mathcal{D}_s^n\left(\frac{\partial L}{\partial \kappa_m}\right)\mathcal{D}_s^{m-1-n}\mathcal{D}_t\kappa\nonumber
\end{eqnarray}
in the boundary terms of (\ref{boundarytermsse2}) show up in the conservation laws (\ref{conslawsse2}), in other words all boundary terms coming from the first set of integration by parts have disappeared. This is no coincidence, it is due to the conflation of $t$ with each group parameter; the proof of this can be found in \cite{Goncalves}.
\end{remark}
 
 In this paper we are interested in showing how the structure of Noether's conservation laws can be used to solve the integration problem for variational problems that are invariant under $SE(2)$ and $SE(3)$. In the present section we have computed Noether's conservation laws for one-dimensional Lagrangians that are invariant under $SE(2)$. In the next section, we will see how these can be used to solve the extremization problems. 

\section{Invariant Lagrangians under $\mathbf{SE(2)}$}\label{SE2}
Recall that the invariantized Euler-Lagrange equation for a one-dimensional Lagrangian invariant under $SE(2)$ is
\begin{equation}\label{se2EL}
\fl\mathsf{E}^y(L)=\mathcal{D}_s^2\mathsf{E}^\kappa(L)+\kappa^2\mathsf{E}^\kappa(L)-\kappa\left(L-\sum_{m=1}\sum_{j=0}^{m-1}(-1)^{j}\mathcal{D}_s^{j}\left(\frac{\partial L}{\partial \kappa_m}\right)\kappa_{m-j}\right),
\end{equation}
and its associated conservation laws are
\begin{eqnarray}
\underbrace{\left(\begin{array}{ccc}
x_s & -y_s & 0\\
y_s & x_s & 0\\
xy_s-yx_s & xx_s+yy_s & 1
\end{array}\right)}_{\mathcal{A}d(\rho(\footnotesize{\zede}))^{-1}}\underbrace{\left(\begin{array}{c}-\lambda(s)-\kappa\mathsf{E}^\kappa(L)\\-\mathcal{D}_s\mathsf{E}^\kappa(L)\\
\mathsf{E}^\kappa(L)\end{array}\right)}_{\boldsymbol{\upsilon}(I)}=\mathbf{c},\nonumber
\end{eqnarray}
where
\begin{eqnarray}
\lambda(s)=-L+\sum_{m=1}\sum_{j=0}^{m-1}(-1)^j\mathcal{D}_s^j\left(\frac{\partial L}{\partial \kappa_m}\right)\kappa_{m-j}.\nonumber
\end{eqnarray}

Multiplying both sides of $\mathcal{A}d(\rho(\zede))^{-1}\boldsymbol{\upsilon}(I)=\mathbf{c}$ by $\mathcal{A}d(\rho(\zede))$ yields the following system of equations
\begin{eqnarray}\label{firsteqse2}
-\lambda(s)-\kappa\mathsf{E}^\kappa(L)=x_sc_1+y_sc_2,\\\label{secondeqse2}
-\mathcal{D}_s\mathsf{E}^\kappa(L)=x_sc_2-y_sc_1,\\\label{thirdeqse2}
\mathsf{E}^\kappa(L)=yc_1-xc_2+c_3.
\end{eqnarray}

Also, we obtain a first integral of the Euler-Lagrange equation as follows. Define
\begin{eqnarray}
\mathsf{B}=\left(\begin{array}{ccc}1 & 0 & 0\\ 0 & 1 & 0\\ 0 & 0 & 0 \end{array}\right),\nonumber
\end{eqnarray}
which satisfies the following equality
\begin{eqnarray}
\mathsf{B}=\mathcal{A}d(\rho)^{-T}\mathsf{B}\mathcal{A}d(\rho)^{-1}.\nonumber
\end{eqnarray}
The first integral of the Euler-Lagrange equation is then
\begin{eqnarray}
\boldsymbol{\upsilon}^T(I)\mathsf{B}\boldsymbol{\upsilon}(I)=\mathbf{c}^T\mathsf{B}\mathbf{c}\nonumber
\end{eqnarray}
i.e.
\begin{equation}\label{se2firstintEL}
(\lambda(s)+\kappa\mathsf{E}^\kappa(L))^2+(\mathcal{D}_s\mathsf{E}^\kappa(L))^2=c_1^2+c_2^2.
\end{equation}
Once $\kappa$ is known, one can see that the integration problem has a straightforward solution. Hence, integrating both sides of Equation (\ref{firsteqse2}) with respect to $s$ yields
\begin{equation}\label{firsteqse2t}
xc_1+yc_2=\int \left[-\lambda(s)-\kappa\mathsf{E}^\kappa(L)\right] \mathrm{d}s,
\end{equation}
and thus solving the two linear equations (\ref{thirdeqse2}) and (\ref{firsteqse2t}) with respect to $x$ and $y$ one obtains
\begin{equation}\label{se2xis}
x(s)=\frac{1}{c_1^2+c_2^2}\left(c_1\int \left[-\lambda(s)-\kappa\mathsf{E}^\kappa(L)\right]\mathrm{d}s-c_2\mathsf{E}^\kappa(L)+c_2c_3\right),
\end{equation}
\begin{equation}\label{se2ipsilon}
y(s)=\frac{1}{c_1^2+c_2^2}\left(c_2\int \left[-\lambda(s)-\kappa\mathsf{E}^\kappa(L)\right]\mathrm{d}s+c_1\mathsf{E}^\kappa(L)+\frac{c_2^2c_3}{c_1}\right).
\end{equation}
One can easily see that Equation (\ref{secondeqse2}) is immediately satisfied as it is the derivative  with respect to $s$ of (\ref{thirdeqse2}).

If we consider the Lagrangian with $L=1$ and plug it in Equations (\ref{se2EL}), (\ref{se2firstintEL}), (\ref{se2xis}) and (\ref{se2ipsilon}), then the solution which minimizes the arc length is the equation of a line, as one would expect. Another famous Lagrangian to consider is
\begin{eqnarray}
\int \kappa^2\mathrm{d}s.\nonumber
\end{eqnarray}
Using the equations above, we obtain, as Euler himself did \cite{Euler}, that the curvature of the minimizing curve satisfies 
\begin{eqnarray}
\kappa_{ss}+\frac{1}{2}\kappa^3=0,\nonumber
\end{eqnarray}
or the first integral of the Euler-Lagrange equation
\begin{eqnarray}
4\kappa_s^2+\kappa^4=c_1^2+c_2^2,\nonumber
\end{eqnarray}
which is solved by an elliptic function. Solutions are known as \emph{Euler's elastica}. For a good historical report see \cite{Levien}.

\section{Invariant Lagrangians under $\mathbf{SE(3)}$}\label{SE3}
In the previous section we showed that for an invariant Lagrangian under $SE(2)$, the invariantized Euler-Lagrange equation and its associated conservation laws could be used to solve the extremising problem. In this section we will proceed analogously and present the solution to one-dimensional variational problems that are invariant under the following $SE(3)$ group action on the $(x(s),y(s),z(s))$-space parametrized by the Euclidean arc length $s$
\begin{equation}\label{gpactse3}
\mathbf{x}\mapsto\widetilde{\mathbf{x}}=\mathsf{R}^{-1}(\mathbf{x}-\mathbf{a}),
\end{equation}
where $\mathsf{R}^{-1}$ represents the rotation in the three-dimensional space
\begin{eqnarray}
\fl\left(\begin{array}{ccc}\cos{\beta}\cos{\gamma}&\cos{\beta}\sin{\gamma} &\sin{\beta} \\
-\sin{\alpha}\sin{\beta}\cos{\gamma}-\cos{\alpha}\sin{\gamma} & -\sin{\alpha}\sin{\beta}\sin{\gamma}+\cos{\alpha}\cos{\gamma}&\sin{\alpha}\cos{\beta} \\
-\cos{\alpha}\sin{\beta}\cos{\gamma}+\sin{\alpha}\sin{\gamma}&-\cos{\alpha}\sin{\beta}\sin{\gamma}-\sin{\alpha}\cos{\gamma} &\cos{\alpha}\cos{\beta}\end{array}\right)\nonumber
\end{eqnarray}
and $\mathbf{a}=(\begin{array}{ccc}a&b&c\end{array})^T$ the translation vector with $\alpha$, $\beta$, $\gamma$, $a$, $b$ and $c$ as the constants that parametrize the group action. 

To solve $SE(3)$ invariant variational problems, we need to find the element $g\in SE(3)$ that sends the tangent to the curve to the $x$-axis, the normal to the curve to the $y$-axis and the point $(x,y,z)$ to the origin, in other words which sends $\zede=(x,y,z,y_s,z_s,z_{ss})$ to the cross section $(0,0,0,0,0,0)$. For that we solve the normalization equations $\widetilde{x}=0$, $\widetilde{y}=0$, $\widetilde{z}=0$, $\widetilde{y_s}=0$, $\widetilde{z_s}=0$ and $\widetilde{z_{ss}}=0$, and thus obtain the right moving frame
\begin{equation}\label{movingse3}
\fl\begin{array}{l}
a=x,\quad b=y,\quad c=z,\quad \alpha=\tan^{-1}\left(\frac{y_s(y_sz_{ss}-z_sy_{ss})-x_s(z_sx_{ss}-x_sz_{ss})}{\sqrt{x_s^2+y_s^2+z_s^2}(x_sy_{ss}-y_sx_{ss})}\right),\\[15pt]
\beta=\tan^{-1}\left(\frac{z_s}{\sqrt{x_s^2+y_s^2}}\right),\quad \gamma=\tan^{-1}\left(\frac{y_s}{x_s}\right).\end{array}
\end{equation}

Consider a one-dimensional variational problem $\mathscr{L}[\mathbf{x}]$ that is invariant under the group action (\ref{gpactse3}). To obtain the invariantized Euler-Lagrange equations, we first rewrite $\mathscr{L}[\mathbf{x}]$ in terms of the generating invariants of the group action, which are the \emph{Euclidean curvature}
\begin{eqnarray}
\kappa=\frac{\parallel\mathbf{x_s}\times\mathbf{x_{ss}}\parallel}{\parallel\mathbf{x_s}\parallel^3},\nonumber
\end{eqnarray}
and \emph{torsion}
\begin{eqnarray}\tau=\frac{\mathbf{x_{sss}}\cdot(\mathbf{x_s}\times\mathbf{x_{ss}})}{\parallel\mathbf{x_s}\times\mathbf{x_{ss}}\parallel^2},\nonumber
\end{eqnarray}
and their derivatives with respect to $s$.

Since $s$ represents the Euclidean arc length, the constraint $\eta=\sqrt{x_s^2+y_s^2+z_s^2}=1$ must be introduced into the variational problem in order to fix parametrization. Hence, the resulting invariantized functional is
\begin{eqnarray}\label{invfunctional2}
\int \left[L(\kappa,\tau,\kappa_s,\tau_s,\kappa_{ss},\tau_{ss},...)-\lambda(s)(\eta-1)\right]\mathrm{d}s,
\end{eqnarray}
where $\lambda(s)$ is a Lagrange multiplier. As for $SE(2)$, this will not reduce the solution set and will simplify the computation of the conservation laws. 

Next we introduce a dummy invariant independent variable $t$ and set $\mathbf{x}=\mathbf{x}(s,t)$ to effect variation. The introduction of a new independent variable results in three new invariants $I^\alpha_t$, for $\alpha=x,y,z$, and a set of syzygies 
\begin{equation}
\mathcal{D}_t\left(\begin{array}{c}\eta\\\kappa\\\tau\end{array}\right)=\mathcal{H}\left(\begin{array}{c}I^x_2\\I^y_2\\I^z_2\end{array}\right),
\end{equation}
where the matrix of operators $\mathcal{H}$ is
\footnotesize{\begin{eqnarray}\fl\left(\kern-6pt\begin{array}{ccc}\mathcal{D}_s & -\kappa & 0\\
\kappa_s & \kappa^2-\tau^2+\mathcal{D}_s^2 & -\tau_s-2\tau\mathcal{D}_s\\
\tau_s+2\tau\mathcal{D}_s & \mathcal{D}_s\left(\frac{\tau_s}{\kappa}\right)+\left(\frac{3\tau_s}{\kappa}-\frac{2\kappa_s\tau}{\kappa^2}\right)\mathcal{D}_s+\frac{2\tau}{\kappa}\mathcal{D}_s^2 & \mathcal{D}_s\left(-\frac{\tau^2}{\kappa}\right)+\left(\kappa-\frac{\tau^2}{\kappa}\right)\mathcal{D}_s-\frac{\kappa_s}{\kappa^2}\mathcal{D}_s^2+\frac{1}{\kappa}\mathcal{D}_s^3\end{array}\kern-6pt\right),\nonumber\end{eqnarray}}\normalsize
where we have already set $\eta=1$.

As in Section \ref{ICV}, we differentiate (\ref{invfunctional2}) with respect to $t$ and then integrate by parts twice to obtain the invariantized Euler-Lagrange equations
\begin{eqnarray}\nonumber
\fl\mathsf{E}^x(L)=&\kappa_s\mathsf{E}^\kappa(L)+\tau_s\mathsf{E}^\tau(L)-\mathcal{D}_s(2\tau\mathsf{E}^\tau(L))+\lambda_s,\nonumber\\[5pt]
\fl\mathsf{E}^y(L)=&\mathcal{D}_s^2\mathsf{E}^\kappa(L)+\frac{2\tau}{\kappa}\mathcal{D}_s^2\mathsf{E}^\tau(L)+\left(\frac{\tau_s}{\kappa}-\frac{2\tau\kappa_s}{\kappa^2}\right)\mathcal{D}_s\mathsf{E}^\tau(L)+(\kappa^2-\tau^2)\mathsf{E}^\kappa(L)\nonumber\\[5pt]
\fl &+2\tau\kappa\mathsf{E}^\tau(L)+\lambda(s)\kappa ,\nonumber\\[5pt]
\fl\mathsf{E}^z(L)=& -\frac{1}{\kappa}\mathcal{D}_s^3\mathsf{E}^\tau(L)+\frac{2\kappa_s}{\kappa^2}\mathcal{D}_s^2\mathsf{E}^\tau(L)+\left(\frac{\kappa_{ss}}{\kappa^2}+\frac{\tau^2}{\kappa}-\frac{2\kappa_s^2}{\kappa^3}-\kappa\right)\mathcal{D}_s\mathsf{E}^\tau(L)\nonumber\\[5pt]
\fl& -\kappa_s\mathsf{E}^\tau(L)+2\tau\mathcal{D}_s\mathsf{E}^\kappa(L)+\tau_s\mathsf{E}^\kappa(L),\nonumber
\end{eqnarray}
and the boundary terms 
\footnotesize{\begin{eqnarray}
\fl\left(2\tau\mathsf{E}^\tau(L)-\lambda\right)I^x_2+\left(\frac{\tau_s}{\kappa}\mathsf{E}^\tau(L)-\frac{2\tau}{\kappa}\mathcal{D}_s\mathsf{E}^\tau(L)-\mathcal{D}_s\mathsf{E}^\kappa(L)\right)I^y_2+\left(\mathsf{E}^\kappa(L)+\frac{2\tau}{\kappa}\mathsf{E}^\tau(L)\right)\mathcal{D}_sI^y_2\nonumber\\
\fl+\left(\kappa\mathsf{E}^\tau(L)-2\tau\mathsf{E}^\kappa(L)-\frac{\tau^2}{\kappa}\mathsf{E}^\tau(L)-\frac{\kappa_s}{\kappa^2}\mathcal{D}_s\mathsf{E}^\tau(L)+\frac{1}{\kappa}\mathcal{D}_s^2\mathsf{E}^\tau(L)\right)I^z_2-\frac{1}{\kappa}\mathcal{D}_s\mathsf{E}^\tau(L)\mathcal{D}_sI^z_2\nonumber\\
\fl+\frac{1}{\kappa}\mathsf{E}^\tau(L)\mathcal{D}_s^2I^z_2=k.\nonumber
\end{eqnarray}}\normalsize

\noindent Note that the boundary terms from the first set of integration by parts disappear after we conflate $t$ with the group parameters.

Using $\mathsf{E}^x(L)=0$ and the fact that

\vspace{-0.4cm}
\footnotesize{\begin{eqnarray}
\fl\kappa_s\mathsf{E}^\kappa(L)+\tau_s\mathsf{E}^\tau(L)=\mathcal{D}_s\Big(L-\sum_{m=1}\sum_{i=0}^{m-1}(-1)^i\mathcal{D}_s^i\frac{\partial L}{\partial \kappa_m}\kappa_{m-i}-\sum_{m=1}\sum_{i=0}^{m-1}(-1)^i\mathcal{D}_s^i\frac{\partial L}{\partial \tau_m}\tau_{m-i}\Big),\nonumber
\end{eqnarray}}\normalsize
(see page 220 of \cite{Mansfield}), we can eliminate $\lambda$. Thus,
\begin{equation}\label{lm2}
\fl\lambda(s)=2\tau\mathsf{E}^\tau(L)-L+\sum_{m=1}\sum_{i=0}^{m-1}(-1)^i\mathcal{D}_s^i\frac{\partial L}{\partial \kappa_m}\kappa_{m-i}+\sum_{m=1}\sum_{i=0}^{m-1}(-1)^i\mathcal{D}_s^i\frac{\partial L}{\partial \tau_m}\tau_{m-i},
\end{equation}
where the constant of integration has been absorbed into the Lagrange multiplier and hence we obtain two Euler-Lagrange equations in two unknowns. 

To calculate the conservation laws associated to the invariantized Euler-Lagrange equations, we must compute the moving frame $\mathcal{A}d(\rho)^{-1}$ and the vector of invariants $\boldsymbol{\upsilon}(I)$. To compute the former we proceed as in Section \ref{NVNcl}: we calculate the Adjoint representation $\mathcal{A}d(g)$ of $SE(3)$ with respect to its generating infinitesimal vector fields
\begin{eqnarray}
\fl\mathbf{v}_a=\partial_x,\;\mathbf{v}_b=\partial_y,\;\mathbf{v}_c=\partial_z,\; \mathbf{v}_\alpha=y\partial_z-z\partial_y,\; \mathbf{v}_\beta=x\partial_z-z\partial_x,\; \mathbf{v}_\gamma=x\partial_y-y\partial_x,\nonumber
\end{eqnarray}
and evaluate it at the frame (\ref{movingse3}), which yields
\begin{eqnarray}
\mathcal{A}d(\rho(\zede))^{-1}=\left(\begin{array}{c|c}
\rho_{FS} & \mathbf{0}\\
\hline
DX\rho_{FS} & D\rho_{FS}D
\end{array}\right),\nonumber
\end{eqnarray}
where
\begin{eqnarray}
\rho_{FS}=\left(\begin{array}{ccc}
\mathbf{x_s} & \frac{\mathbf{x_{ss}}}{\kappa} & \frac{\mathbf{x_s}\times\mathbf{x_{ss}}}{\kappa}\end{array}\right)\nonumber
\end{eqnarray}
is the Frenet-Serret frame, $D$ is the following diagonal matrix
\begin{eqnarray}
D=\left(\begin{array}{ccc}
1&0&0\\
0&-1&0\\
0&0&1\end{array}\right),\nonumber
\end{eqnarray}
and $X$ is the matrix
\begin{eqnarray}
X=\left(\begin{array}{ccc}
0 & -z & y\\
z & 0 & -x\\
-y & x & 0\end{array}\right).\nonumber
\end{eqnarray}

As seen in Section \ref{NVNcl}, to obtain the vector of invariants we must first find the boundary terms that are linear in the $I^\alpha_{2K}$. To do so, consider the boundary terms obtained from the calculation of the invariantized Euler-Lagrange equations

\vspace{-0.4cm}
\footnotesize{\begin{eqnarray}\label{bts}
\fl\begin{array}{l}\left(2\tau\mathsf{E}^\tau(L)-\lambda\right)I^x_2+\left(\frac{\tau_s}{\kappa}\mathsf{E}^\tau(L)-\frac{2\tau}{\kappa}\mathcal{D}_s\mathsf{E}^\tau(L)-\mathcal{D}_s\mathsf{E}^\kappa(L)\right)I^y_2+\left(\mathsf{E}^\kappa(L)+\frac{2\tau}{\kappa}\mathsf{E}^\tau(L)\right)\mathcal{D}_sI^y_2\\
+\left(\kappa\mathsf{E}^\tau(L)-2\tau\mathsf{E}^\kappa(L)-\frac{\tau^2}{\kappa}\mathsf{E}^\tau(L)-\frac{\kappa_s}{\kappa^2}\mathcal{D}_s\mathsf{E}^\tau(L)+\frac{1}{\kappa}\mathcal{D}_s^2\mathsf{E}^\tau(L)\right)I^z_2-\frac{1}{\kappa}\mathcal{D}_s\mathsf{E}^\tau(L)\mathcal{D}_sI^z_2\\
+\frac{1}{\kappa}\mathsf{E}^\tau(L)\mathcal{D}_s^2I^z_2=k.\end{array}
\end{eqnarray}}\normalsize
Substituting $\mathcal{D}_sI^y_2$, $\mathcal{D}_sI^z_2$ and $\mathcal{D}_s^2I^z_2$ in (\ref{bts}) by the differential formulae
\begin{eqnarray}
\mathcal{D}_s I_2^y=-\kappa I_2^x+I_{12}^y+\tau I_2^z,\nonumber\\
\mathcal{D}_s I_2^z=-\tau I_2^y+I_{12}^z,\nonumber\\
\mathcal{D}_s^2 I_2^z=\tau\kappa I_2^x-\tau_s I_2^y-2\tau I_{12}^y-\tau^2 I_2^z+I_{112}^z,\nonumber
\end{eqnarray}
which were obtained from (\ref{diffinvariant}), yields that the boundary terms are linear in the $I^\alpha_{2K}$,
\footnotesize{\begin{eqnarray}
\fl\begin{array}{l}
\left(\begin{array}{c} I^x_2\end{array}\right)
\underbrace{\left(\begin{array}{c}\displaystyle{-\kappa\mathsf{E}^\kappa(L)+\tau\mathsf{E}^\tau(L)-\lambda(s)}\end{array}\right)}_{\mathcal{C}^x}+\left(\begin{array}{cc}I^y_2 & I^y_{12} \end{array}\right)\underbrace{\left(\begin{array}{c} \displaystyle{-\mathcal{D}_s\mathsf{E}^\kappa(L)-\frac{\tau}{\kappa}\mathcal{D}_s\mathsf{E}^\tau(L)}\\[5pt]\mathsf{E}^\kappa(L)\end{array}\right)}_{\mathcal{C}^y}\\[15pt]
+\left(\begin{array}{ccc}I^z_2 & I^z_{12} & I^z_{112}\end{array}\right)\underbrace{\left(\begin{array}{c} \displaystyle{\frac{1}{\kappa}\mathcal{D}^2_s\mathsf{E}^\tau(L)-\frac{\kappa_s}{\kappa^2}\mathcal{D}_s\mathsf{E}^\tau(L)+\kappa\mathsf{E}^\tau(L)-\tau\mathsf{E}^\kappa(L)}\\[5pt]\displaystyle{-\frac{1}{\kappa}\mathcal{D}_s\mathsf{E}^\tau(L)}\\[5pt]\displaystyle{\frac{1}{\kappa}\mathsf{E}^\tau(L)}\end{array}\right)}_{\mathcal{C}^z}=k.\nonumber
\end{array}
\end{eqnarray}}\normalsize 
Finally, adding the products of the matrices of invariantized infinitesimals $\Omega^\alpha(I)$ with the vectors $\mathcal{C}^\alpha$ yields the vector of invariants
\begin{eqnarray}
\boldsymbol{\upsilon}(I)=\left(\begin{array}{c}
\tau\mathsf{E}^\tau(L)-\kappa\mathsf{E}^\kappa(L)-\lambda(s)\\
-\mathcal{D}_s\mathsf{E}^\kappa(L)-\frac{\tau}{\kappa}\mathcal{D}_s\mathsf{E}^\tau(L)\\
\frac{1}{\kappa}\mathcal{D}_s^2\mathsf{E}^\tau(L)-\frac{\kappa_s}{\kappa^2}\mathcal{D}_s\mathsf{E}^\tau(L)+\kappa\mathsf{E}^\tau(L)-\tau\mathsf{E}^\kappa(L)\\
\mathsf{E}^\tau(L)\\
-\frac{1}{\kappa}\mathcal{D}_s\mathsf{E}^\tau(L)\\
\mathsf{E}^\kappa(L)
\end{array}\right),\nonumber
\end{eqnarray}
where the $\Omega^\alpha(I)$ are
\begin{eqnarray}
\Omega^x(I)=\left(\begin{array}{c}1\\0\\0\\0\\0\\0\end{array}\right),\quad \Omega^y(I)=\left(\begin{array}{cc}0 & 0\\1 & 0\\0 & 0\\0 & 0\\0 & 0\\0 & 1\end{array}\right),\quad\Omega^z(I)=\left(\begin{array}{ccc}0 & 0 & 0\\0 & 0 & 0\\1 & 0 & 0\\0 & 0 & \kappa\\0 & 1 & 0\\0 & 0 & 0\end{array}\right).\nonumber
\end{eqnarray}

Thus, the conservation laws are
\begin{equation}\label{consse3}
\fl\footnotesize{\left(\begin{array}{c|c}
\rho_{FS} & \mathbf{0}\\
\hline
DX\rho_{FS} & D\rho_{FS}D
\end{array}\right)\left(\begin{array}{c}
\tau\mathsf{E}^\tau(L)-\kappa\mathsf{E}^\kappa(L)-\lambda(s)\\
-\mathcal{D}_s\mathsf{E}^\kappa(L)-\frac{\tau}{\kappa}\mathcal{D}_s\mathsf{E}^\tau(L)\\
\frac{1}{\kappa}\mathcal{D}_s^2\mathsf{E}^\tau(L)-\frac{\kappa_s}{\kappa^2}\mathcal{D}_s\mathsf{E}^\tau(L)+\kappa\mathsf{E}^\tau(L)-\tau\mathsf{E}^\kappa(L)\\
\mathsf{E}^\tau(L)\\
-\frac{1}{\kappa}\mathcal{D}_s\mathsf{E}^\tau(L)\\
\mathsf{E}^\kappa(L)
\end{array}\right)=\left(\begin{array}{c}
\mathbf{c_1}\\
\mathbf{c_2}
\end{array}\right)=\mathbf{c},}
\end{equation}\normalsize
where $\mathbf{c_1}=(c_1,c_2,c_3)^T$, $\mathbf{c_2}=(c_4,c_5,c_6)^T$ are constant vectors and $\lambda(s)$ is equal to (\ref{lm2}).

We shall see in the remainder of this section how the conservation laws (\ref{consse3}) can help reduce the integration problem.

To demonstrate this, we will start by simplifying the conservation laws (\ref{consse3}) in two steps. These simplifications will then lead to an overdetermined system of equations for $x$, $y$ and $z$, which will be solved with relative ease once $\kappa$ and $\tau$ are known. Finally, we will give a reason behind the choice in the order in which we solve the equations.

In the first step of the simplification, we apply an element of $SE(3)$, say $\mathcal{A}d(g)^{-1}$, to both sides of $\mathcal{A}d(\rho(\zede))^{-1}\boldsymbol{\upsilon}(I)=\mathbf{c}$ such that it maps $\mathbf{c}_1$ and $\mathbf{c}_2$ to the $z$-axis. Let $\mathcal{A}d(g)$ act on $\mathbf{c}$ as follows
\begin{eqnarray}
\mathcal{A}d(g)\mathbf{c}=\left(\begin{array}{cc}\mathsf{R} & \mathbf{0}\\
DA\mathsf{R} & D\mathsf{R}D\end{array}\right)\left(\begin{array}{c}\mathbf{c}_1\\\mathbf{c}_2\end{array}\right)=\left(\begin{array}{c}\widetilde{\mathbf{c}_1}\\\widetilde{\mathbf{c}_2}\end{array}\right),\nonumber
\end{eqnarray}
where $\mathsf{R}$ is the three-dimensional rotation
\begin{eqnarray}
\fl\mathsf{R}=\left(\begin{array}{ccc}\cos{\beta}\cos{\gamma} & -\sin{\alpha}\sin{\beta}\cos{\gamma}-\cos{\alpha}\sin{\gamma} & -\cos{\alpha}\sin{\beta}\cos{\gamma}+\sin{\alpha}\sin{\gamma}\\
\cos{\beta}\sin{\gamma} & -\sin{\alpha}\sin{\beta}\sin{\gamma}+\cos{\alpha}\cos{\gamma} & -\cos{\alpha}\sin{\beta}\sin{\gamma}-\sin{\alpha}\cos{\gamma}\\
\sin{\beta} & \sin{\alpha}\cos{\beta} & \cos{\alpha}\cos{\beta}\end{array}\right),\nonumber
\end{eqnarray}
$D$ is the diagonal matrix, $D=\mathrm{diag}(1,-1,1)$, and $A$ the matrix, 
\begin{eqnarray}
A=\left(\begin{array}{ccc}0& -c & b\\c & 0 & -a\\-b & a & 0\end{array}\right).\nonumber
\end{eqnarray}
We can easily verify that the Adjoint representation of $SE(3)$ does not act freely on the constant vector $\mathbf{c}$, since it preserves the length of $\mathbf{c}_1$ and the quantity $\mathbf{c}_1^TD\mathbf{c}_2$. Indeed to prove the latter, we multiply through
\begin{eqnarray}
DA\mathsf{R}\mathbf{c}_1+D\mathsf{R}D\mathbf{c}_1=\widetilde{\mathbf{c}_2}\nonumber
\end{eqnarray}
by $\mathbf{c}_1^T\mathsf{R}^TD$, and then obtain that
\begin{eqnarray}
\underbrace{\widetilde{\mathbf{c}_1}^TA\widetilde{\mathbf{c}_1}}_{=0}+\mathbf{c}_1^TD\mathbf{c}_2=\widetilde{\mathbf{c}_1}^TD\widetilde{\mathbf{c}_2}.\nonumber
\end{eqnarray}
Thus, let $\mathcal{A}d(g)^{-1}$ send $\mathbf{c}$ to
\begin{eqnarray}
\mathbf{C}=\left(\begin{array}{cccccc} 0 & 0 & |\mathbf{c}_1| & 0 & 0 & \frac{\mathbf{c}_1^TD\mathbf{c}_2}{|\mathbf{c}_1|}\end{array}\right)^T.\nonumber
\end{eqnarray}
Applying $\mathcal{A}d(g)^{-1}$ to the conservation laws yields
\begin{eqnarray}
\mathcal{A}d(g)^{-1}\mathcal{A}d(\rho(\zede))^{-1}\boldsymbol{\upsilon}(I)=\mathcal{A}d(g)^{-1}\mathbf{c},\nonumber
\end{eqnarray}
which reduces to
\begin{equation}\label{interm1}
\mathcal{A}d(\rho(\widetilde{\zede}))^{-1}\boldsymbol{\upsilon}(I)=\mathbf{C}
\end{equation}
by the equivariance of the right moving frame $\mathcal{A}d(\rho(\zede))^{-1}$.

The second step of our simplification consists of applying $\mathcal{A}d(\rho(\widetilde{\zede}))$ to (\ref{interm1}) to obtain the following system of equations
\begin{equation}\label{IN1}
\fl |\mathbf{c_1}|\widetilde{z_s}=\upsilon^{(1)}(I),
\end{equation}
\begin{equation}\label{IN2}
\fl \frac{|\mathbf{c_1}|}{\kappa}\widetilde{z_{ss}}=\upsilon^{(2)}(I),
\end{equation}
\begin{equation}\label{IN3}
\fl \frac{|\mathbf{c_1}|}{\kappa}(\widetilde{x_s}\widetilde{y_{ss}}-\widetilde{y_s}\widetilde{x_{ss}})=\upsilon^{(3)}(I),
\end{equation}
\begin{equation}\label{IN4}
\fl |\mathbf{c_1}|(\widetilde{x}\widetilde{y_s}-\widetilde{y}\widetilde{x_s})+\frac{\mathbf{c_1}^TD\mathbf{c_2}}{|\mathbf{c_1}|}\widetilde{z_s}=\upsilon^{(4)}(I),
\end{equation}
\begin{equation}\label{IN5}
\fl \frac{|\mathbf{c_1}|}{\kappa}(\widetilde{x_{ss}}\widetilde{y}-\widetilde{y_{ss}}\widetilde{x})-\frac{\mathbf{c_1}^TD\mathbf{c_2}}{\kappa|\mathbf{c_1}|}\widetilde{z_{ss}}=\upsilon^{(5)}(I),
\end{equation}
\begin{equation}\label{IN6}
\fl \frac{|\mathbf{c_1}|}{\kappa}(\widetilde{x}(\widetilde{z_s}\widetilde{x_{ss}}-\widetilde{x_s}\widetilde{z_{ss}})-\widetilde{y}(\widetilde{y_s}\widetilde{z_{ss}}-\widetilde{z_s}\widetilde{y_{ss}}))+\frac{\mathbf{c_1}^TD\mathbf{c_2}}{\kappa|\mathbf{c_1}|}(\widetilde{x_s}\widetilde{y_{ss}}-\widetilde{y_s}\widetilde{x_{ss}})=\upsilon^{(6)}(I),
\end{equation}
where we have used $\upsilon^{(j)}(I)$ to denote the $j$-th component of $\boldsymbol{\upsilon}(I)$. This overdetermined system of equations can now be solved more easily.

The two first integrals of the Euler-Lagrange equations are obtained as follows. Define $\mathsf{B}=\mathrm{diag}(\begin{array}{cccccc}1&1&1&0&0&0\end{array})$, which satisfies $\mathsf{B}=\mathcal{A}d(\rho)^{-T}\mathsf{B}\mathcal{A}d(\rho)^{-1}$. Then we obtain the first integral of the Euler-Lagrange equations,
$$\boldsymbol{\upsilon}^T(I)\mathsf{B}\boldsymbol{\upsilon}(I)=\mathbf{C}^T\mathsf{B}\mathbf{C},$$
which is equaivalent to
\begin{eqnarray}\label{1firstintegral}
\scriptstyle{\fl\begin{array}{l}
\left(\tau\mathsf{E}^\tau(L)-\kappa\mathsf{E}^\kappa(L)-\lambda(s)\right)^2+\left(-\mathcal{D}_s\mathsf{E}^\kappa(L)-\frac{\tau}{\kappa}\mathcal{D}_s\mathsf{E}^\tau(L)\right)^2\\[10pt]
+\left(\frac{1}{\kappa}\mathcal{D}_s^2\mathsf{E}^\tau(L)-\frac{\kappa_s}{\kappa^2}\mathcal{D}_s\mathsf{E}^\tau(L)+\kappa\mathsf{E}^\tau(L)-\tau\mathsf{E}^\kappa(L)\right)^2=c_1^2+c_2^2+c_3^2.\end{array}}
\end{eqnarray}
For the second first integral, we define
$$\mathsf{D}=\left(\begin{array}{cc}\mathbf{0} & D\\ D& \mathbf{0}\end{array}\right),\qquad D=\mathrm{diag}(\begin{array}{ccc}1 &-1 &1\end{array}),$$
which satifies $\mathsf{D}=\mathcal{A}d(\rho)^{-T}\mathsf{D}\mathcal{A}d(\rho)^{-1}$. The first integral of the Euler-Lagrange equations is then
$$\boldsymbol{\upsilon}^T(I)\mathsf{D}\boldsymbol{\upsilon}(I)=\mathbf{C}^T\mathsf{D}\mathbf{C},$$
i.e.
\begin{eqnarray}\label{2firstintegral}
\scriptstyle{\fl \begin{array}{l}
\left(\tau\mathsf{E}^\tau(L)-\kappa\mathsf{E}^\kappa(L)-\lambda(s)\right)\mathsf{E}^\tau(L)+\left(-\mathcal{D}_s\mathsf{E}^\kappa(L)-\frac{\tau}{\kappa}\mathcal{D}_s\mathsf{E}^\tau(L)\right)\frac{1}{\kappa}\mathcal{D}_s\mathsf{E}^\tau(L)\\[10pt]
+\left(\frac{1}{\kappa}\mathcal{D}_s^2\mathsf{E}^\tau(L)-\frac{\kappa_s}{\kappa^2}\mathcal{D}_s\mathsf{E}^\tau(L)+\kappa\mathsf{E}^\tau(L)-\tau\mathsf{E}^\kappa(L)\right)\mathsf{E}^\kappa(L)=c_1c_4-c_2c_5+c_3c_6.
\end{array}}
\end{eqnarray}

We can use the above first integrals to determine $\kappa$ and $\tau$. Once we have solved for these, we use the system of simplified conservation laws to solve for $\widetilde{x}$, $\widetilde{y}$ and $\widetilde{z}$. 

Hence, solving Equation (\ref{IN1}) gives
\begin{eqnarray}\label{zvariable}
\widetilde{z(s)}=\frac{1}{|\mathbf{c_1}|}\displaystyle{\int \upsilon^{(1)}(I)\mathrm{d}s}.
\end{eqnarray}
Next, multiplying Equation (\ref{IN3}) by $-\frac{\mathbf{c_1}^TD\mathbf{c_2}}{|\mathbf{c_1}|^2}$ and adding it to Equation (\ref{IN6}) yields
\begin{eqnarray}
\fl \frac{|\mathbf{c_1}|}{\kappa}(\widetilde{x}(\widetilde{z_s}\widetilde{x_{ss}}-\widetilde{x_s}\widetilde{z_{ss}})-\widetilde{y}(\widetilde{y_s}\widetilde{z_{ss}}-\widetilde{z_s}\widetilde{y_{ss}}))=\upsilon^{(6)}(I)-\frac{\mathbf{c_1}^TD\mathbf{c_2}}{|\mathbf{c_1}|^2}\upsilon^{(3)}(I),\nonumber
\end{eqnarray}
which simplifies to
\begin{eqnarray}\nonumber
\fl|\mathbf{c_1}|\left(\widetilde{z_s}\left(\frac{1}{2}\mathcal{D}_s^2(\widetilde{\mathbf{x}}\cdot\widetilde{\mathbf{x}})-1\right)-\frac{1}{2}\widetilde{z_{ss}}\mathcal{D}_s(\widetilde{\mathbf{x}}\cdot\widetilde{\mathbf{x}})\right)=\kappa\upsilon^{(6)}(I)-\kappa\frac{\mathbf{c_1}^TD\mathbf{c_2}}{|\mathbf{c_1}|^2}\upsilon^{(3)}(I),\end{eqnarray}
where $\frac{1}{2}\mathcal{D}_s^2(\widetilde{\mathbf{x}}\cdot\widetilde{\mathbf{x}})-1=\widetilde{\mathbf{x}}\cdot\widetilde{\mathbf{x_{ss}}}$ and $\frac{1}{2}\mathcal{D}_s(\widetilde{\mathbf{x}}\cdot\widetilde{\mathbf{x}})=\widetilde{\mathbf{x}}\cdot\widetilde{\mathbf{x_s}}$. Setting $\mathcal{D}_s(\widetilde{\mathbf{x}}\cdot\widetilde{\mathbf{x}})=h(s)$ and substituting $\widetilde{z_s}$ by (\ref{IN1}) and $\widetilde{z_{ss}}$ by its derivative yields a linear equation for $h$
\begin{eqnarray}\nonumber
\fl\mathcal{D}_s h-\frac{\mathcal{D}_s\upsilon^{(1)}(I)}{\upsilon^{(1)}(I)}h=\left.2\kappa\left(\upsilon^{(6)}(I)-\frac{\mathbf{c_1}^TD\mathbf{c_2}}{|\mathbf{c_1}|^2}\upsilon^{(3)}(I)\right)\middle/\upsilon^{(1)}(I)\right.+2.
\end{eqnarray}
Solving for $h$ we obtain
\begin{eqnarray}\nonumber
\fl h(s)=\upsilon^{(1)}(I)\int\frac{1}{\upsilon^{(1)}(I)}\left(\left.2\kappa\left(\upsilon^{(6)}(I)-\frac{\mathbf{c_1}^TD\mathbf{c_2}}{|\mathbf{c_1}|^2}\upsilon^{(3)}(I)\right)\middle/\upsilon^{(1)}(I)\right.+2\right)\mathrm{d}s.
\end{eqnarray}
Hence,
\begin{equation}\label{innerprod}
\fl\mathcal{D}_s(\widetilde{\mathbf{x}}\cdot\widetilde{\mathbf{x}})=\upsilon^{(1)}(I)\int\frac{1}{\upsilon^{(1)}(I)}\left(\left.2\kappa\left(\upsilon^{(6)}(I)-\frac{\mathbf{c_1}^TD\mathbf{c_2}}{|\mathbf{c_1}|^2}\upsilon^{(3)}(I)\right)\middle/\upsilon^{(1)}(I)\right.+2\right)\mathrm{d}s.
\end{equation}
To solve Equations (\ref{IN4}) and (\ref{innerprod}), we use the cylindrical coordinates
\begin{eqnarray}
\widetilde{x(s)}=r(s)\cos{\theta(s)},\quad \widetilde{y(s)}=r(s)\sin{\theta(s)},\quad \widetilde{z(s)}=\widetilde{z(s)}.\nonumber
\end{eqnarray}
Starting with Equation (\ref{innerprod}), we obtain
\begin{eqnarray} \label{radius} r(s)^2=\int h(s)\mathrm{d}s-\widetilde{z(s)}^2,\end{eqnarray}
as $\mathcal{D}_s(\widetilde{\mathbf{x}}\cdot\widetilde{\mathbf{x}})=\mathcal{D}_s\Big(r(s)^2+\widetilde{z(s)}^2\Big)$. After applying the change of coordinates, Equation (\ref{IN4}) becomes
\begin{eqnarray}\nonumber
r(s)^2\theta_s=\frac{1}{|\mathbf{c_1}|}\left(\upsilon^{(4)}(I)-\frac{\mathbf{c_1}^TD\mathbf{c_2}}{|\mathbf{c_1}|^2}\upsilon^{(1)}(I)\right),
\end{eqnarray}
 and then solving for $\theta$ yields
\begin{eqnarray}\label{angle}
\theta(s)=\int\frac{1}{r(s)^2|\mathbf{c_1}|}\left(\upsilon^{(4)}(I)-\frac{\mathbf{c_1}^TD\mathbf{c_2}}{|\mathbf{c_1}|^2}\upsilon^{(1)}(I)\right)\mathrm{d}s.
\end{eqnarray}
To recover $x$, $y$ and $z$, we act on $\widetilde{x}$, $\widetilde{y}$ and $\widetilde{z}$ as follows
$$\widetilde{\mathbf{x}}\mapsto\mathbf{x}=\mathsf{R}\widetilde{\mathbf{x}}+\mathbf{a},$$
where $\mathsf{R}$ is the three-dimensional rotation
\begin{eqnarray}
\fl\mathsf{R}\kern-2pt=\kern-2pt\left(\begin{array}{ccc}\cos{\beta}\cos{\gamma} & -\sin{\alpha}\sin{\beta}\cos{\gamma}-\cos{\alpha}\sin{\gamma} & -\cos{\alpha}\sin{\beta}\cos{\gamma}+\sin{\alpha}\sin{\gamma}\\
\cos{\beta}\sin{\gamma} & -\sin{\alpha}\sin{\beta}\sin{\gamma}+\cos{\alpha}\cos{\gamma} & -\cos{\alpha}\sin{\beta}\sin{\gamma}-\sin{\alpha}\cos{\gamma}\\
\sin{\beta} & \sin{\alpha}\cos{\beta} & \cos{\alpha}\cos{\beta}\end{array}\right),\nonumber
\end{eqnarray}
and $\mathbf{a}=(\begin{array}{ccc} a & b & c\end{array})^T$ is the translation vector, with
\begin{eqnarray}
\fl\alpha=-\tan^{-1}\left(\frac{\sqrt{|\mathbf{c_1}|^2\cos^2{\beta}-c_3^2}}{c_3}\right),\; \gamma=\tan^{-1}\left(\frac{c_2c_3\sin{\beta}+c_1\sqrt{|\mathbf{c_1}|^2\cos^2{\beta}-c_3^2}}{c_1c_3\sin{\beta}-c_2\sqrt{|\mathbf{c_1}|^2\cos^2{\beta}-c_3^2}}\right),\nonumber
\end{eqnarray}
\begin{eqnarray}
\fl a=\frac{c_1}{c_3}c+\frac{c_5|\mathbf{c_1}|^2+c_2\mathbf{c_1}^TD\mathbf{c_2}}{c_3|\mathbf{c_1}|^2},\quad b=\frac{c_2}{c_3}c+\frac{c_4|\mathbf{c_1}|^2-c_1\mathbf{c_1}^TD\mathbf{c_2}}{c_3|\mathbf{c_1}|^2},\nonumber
\end{eqnarray}
and where $\beta$ and $c$ are free.

Although only four of the equations of the system were used to solve for $x$, $y$ and $z$, we know that the remaining two equations have been satisfied. Indeed, if we differentiate $\mathcal{A}d(\rho(\widetilde{\zede}))^{-1}\boldsymbol{\upsilon}(I)=\mathbf{C}$ with respect to $s$  and multiply by $\mathcal{A}d(\rho(\widetilde{\zede}))$, then we get
\begin{eqnarray}
\mathcal{D}_s\boldsymbol{\upsilon}(I)=\mathcal{D}_s\left(\mathcal{A}d(\rho(\widetilde{\zede}))\right))\mathcal{A}d(\rho(\widetilde{\zede}))^{-1}\boldsymbol{\upsilon}(I),\nonumber
\end{eqnarray}
which is equivalent to
\begin{equation}\label{eliminationideal}
\mathcal{D}_s\boldsymbol{\upsilon}(I)=\left(\begin{array}{cccccc}
0 & \kappa & 0 & 0 & 0 & 0\\
-\kappa & 0 & \tau & 0 & 0 &  0\\
0 & -\tau & 0 & 0 & 0 & 0\\
0 & 0 & 0 & 0 & -\kappa & 0\\
0 & 0 & -1 & \kappa & 0 &-\tau\\
0 & -1 & 0 & 0 & \tau & 0
\end{array}\right)\boldsymbol{\upsilon}(I).
\end{equation}
The above system of equations not only forms part of an \emph{elimination ideal} as it only involves invariants, but because the invariants appear only on the right-hand sides of the Equations (\ref{IN1}), (\ref{IN2}), (\ref{IN3}), (\ref{IN4}), (\ref{IN5}) and (\ref{IN6}), they also encode the relationships between the equations themselves. Thus we see ${\cal D}_s (\mathrm{Equation}\;(\ref{IN1}))=\kappa (\mathrm{Equation}\;(\ref{IN2}))$ and so forth. Using the Equations in (\ref{eliminationideal}) we can eliminate (\ref{IN2}) and (\ref{IN5}) from the system.

In the following two well-known examples, we can verify that we obtain the expected solutions.

\begin{example}
For the Lagrangian with $L=1$, Equations (\ref{zvariable}), (\ref{radius}) and (\ref{angle}) yield the equation of a line in parametric form as the solution that minimizes the arc length, as expected. Note that $|\widetilde{\mathbf{x}_s}|=1$ imposes a condition on the constants of integration.
\end{example}

\begin{example}
Consider the Lagrangian $\int \kappa^2\mathrm{d}s$ with torsion, $\tau$, equal to zero. Then we obtain that $\kappa$ satisfies the Euler-Lagrange equation
\begin{eqnarray}
\kappa_{ss}+\frac{1}{2}\kappa^3=0,\nonumber
\end{eqnarray}
which is the same equation as for the $SE(2)$ case, or the first integral 
\begin{eqnarray}
\kappa^4+4\kappa_s^2=|\mathbf{c}_1|.\nonumber
\end{eqnarray}
From Equation (\ref{2firstintegral}) we know that
\begin{eqnarray}
\mathbf{c}_1^TD\mathbf{c}_2=0,\nonumber
\end{eqnarray}
and from Equations (\ref{zvariable}), (\ref{radius}) and (\ref{angle}) we obtain
\begin{eqnarray}
\widetilde{z(s)}=-\frac{1}{|\mathbf{c}_1|}\int\kappa^2\mathrm{d}s,\nonumber\\[7pt]
r(s)^2=-\int\Big[\kappa^2\int\frac{2}{\kappa^2}\mathrm{d}s\Big]\mathrm{d}s-\frac{1}{|\mathbf{c}_1|^2}\left(\int\kappa^2\mathrm{d}s\right)^2,\nonumber\\[8pt]
\theta(s)=A,\nonumber
\end{eqnarray}
where $A$ is a constant. So the solution $\widetilde{\mathbf{x}}$ lies on a plane that includes the $\widetilde{z}$-axis, as expected.
\end{example}

\section{Conclusion}
Noether's First Theorem is a well-known result which provides conservation laws for Lie group invariant variational problems. In recent work \cite{GoncalvesMansfield}, the mathematical structure of both the Euler-Lagrange system and the set of conservation laws was given in terms of the differential invariants of the group action and a moving frame. It is the knowledge of this structure that allows one to solve the invariant variational problems under some group action with relative ease. In this paper, we examine one-dimensional variational problems that are invariant under the actions of $SE(2)$ and $SE(3)$. For both cases, we obtain the invariantized Euler-Lagrange equations and their associated conservation laws in the new format, from which we can then obtain the solution to the variational problem by quadratures. One can verify that this method leads to a far simpler computational problem than the one given in the original variables.

\pagebreak
\noindent \textbf{References}
\vspace{0.4cm}
\bibliographystyle{unsrt}
\bibliography{jpA}
\end{document}